\documentclass[12pt,twoside]{article}

\date{}

\begin{document}

\centerline{\bf International Journal of Contemporary
Mathematical Sciences, Vol. x, 200x, no. xx, xxx - xxx}
\centerline{}
 \centerline{}
 \centerline{\Large{\bf On the Order of Schur
Multipliers of Finite Abelian ${\bf p}$-Groups}}
 \centerline{}

\centerline{\bf {B. Mashayekhy*, F. Mohammadzadeh, and A.
Hokmabadi}}

\centerline{}

\centerline{Department of Mathematics,} \centerline{ Center of
Excellence in Analysis on Algebraic Structures,} \centerline{
Ferdowsi University of Mashhad,} \centerline{ P. O. Box
1159-91775, Mashhad, Iran.}\centerline{* Correspondence:
mashaf@math.um.ac.ir}

\newtheorem{Theorem}{\quad Theorem}[section]

\newtheorem{Definition}[Theorem]{\quad Definition}

\newtheorem{Corollary}[Theorem]{\quad Corollary}

\newtheorem{Lemma}[Theorem]{\quad Lemma}

\newtheorem{Example}[Theorem]{\quad Example}

\begin{abstract}
Let $G$ be a finite $p$-group of order $p^{n}$ with
$|M(G)|=p^{\frac{n(n-1)}{2}-t},$ where $M(G)$ is the Schur
multiplier of $G$. Ya.G. Berkovich, X. Zhou, and G. Ellis have
determined the structure of $G$ when $t=0,1,2,3$. In this paper,
we are going to find some structures for an abelian $p$-group $G$
with conditions on the exponents of $G, M(G),$ and $S_2M(G)$,
where $S_2M(G)$ is the metabelian multiplier of $G$.
\end{abstract}



\section{Introduction and Preliminaries}

 Let $G$ be any group with a presentation $G\cong F/R$, where F is
a free group. Then the Baer invariant of $G$ with respect to the
variety of groups $\mathcal{V}$, denoted by ${\mathcal V}M(G)$,
is defined to be
$${\mathcal V}M(G)=\frac{R \cap V(F)}{[RV^{*}F]},$$
where $V$ is the set of words of the variety $\mathcal{V}$, $V(F)$
is the verbal subgroup of $F$ and
\newpage
$$ [R V^{*}F]=\langle v(f_1,...,f_{i-1},f_{i}r,f_{i+1},...,f_n)
 v(f_1,...,f_{i},...,f_n)^{-1} |$$ $$ r\in R, f_i\in F, v
\in V, 1\leq i\leq n, n \in N \rangle.$$ In particular, if
$\mathcal{V}$ is the variety of abelian groups, $\mathcal{A}$,
then the Baer invariant of the group $G$ will be $(R \cap F')/[R,
F]$ which is isomorphic to the well-known notion the Schur
multiplier of $G$, denoted by $M(G)$
(see [5,6] for further details).\\
If $\mathcal{V}$ is the variety of polynilpotent groups of class
row $(c_1,...,c_t)$, $ {\mathcal N}_{c_1,c_2,...,c_t}$, then the
Baer invariant of a group $G$ with respect to this variety is as
follows:
\begin{equation} {\mathcal N}_{c_1,c_2,...,c_t} M(G)=\frac{R \cap \gamma_{c_t+1}\circ
...\circ \gamma_{c_1+1}(F)}{[R,\ _{c_1}F,\
_{c_2}\gamma_{c_1+1}(F),...,\ _{c_t}\gamma_{c_{t-1}+1}\circ
...\circ \gamma_{c_1+1}(F)]},
\end{equation} where
$\gamma_{c_t+1}\circ ...\circ
\gamma_{c_1+1}(F)=\gamma_{c_t+1}(\gamma_{c_{t-1}+1}( ...
(\gamma_{c_1+1}(F))...))$ are the term of iterated lower central
series of $F$. See [4] for the equality
$$[R{\mathcal{N}}^{*}_{c_1,...,c_t}F]=[R,\ _{c_1}F,\
_{c_2}\gamma_{c_1+1}(F),...,\ _{c_t}\gamma_{c_{t-1}+1}\circ
...\circ \gamma_{c_1+1}(F)].$$ In particular, if $c_i=1$ for $1
\leq i \leq t$, then ${\mathcal N}_{c_1,c_2,...,c_t}$ is the
variety of
solvable groups of length at most $t \geq 1, {\mathcal S}_t.$\\
In 1956, J.A. Green [3] showed that the order of the Schur
multiplier of a finite $p$-group of order $p^{n}$ is bounded by
$p^{\frac{n(n-1)}{2}}$, and hence equals to
$p^{\frac{n(n-1)}{2}-t}$, for some nonnegative integer $t$. In
1991, Ya.G. Berkovich [1] has determined all finite $p$-groups
$G$ for which $t=0,1$. The groups for which $t=0$ are exactly
elementary ablian $p$-groups, and the groups for which $t=1$ are
cyclic groups of order $p^{2}$ or the nonabelian group of order
$p^{3}$ with exponent $p>2$. In 1994, X. Zhou [7] found all finite
$p$-groups for $t=2$. He showed that these groups are the direct
product of two cyclic groups of order $p^{2}$ and $p$ or the
direct product of a cyclic group of order $p$ and the nonabelian
group of order $p^{3}$ and exponent $p>2$ or the dihedral group
of order 8. G. Ellis [2] determined all finite $p$-groups $G$
with $t=0,1,2,3$ in a quite different method to that of [1] and
[7] as follows:
\begin{Theorem} ([2]). Let $G$ be a group of
prime-power order $p^{n}$. Suppose that $M(G)$ has order
$p^{\frac{n(n-1)}{2}-t}$. Then $t \geq 0$ and \\
(i) $t=0$ if and only if $G$ is elementary abelian;\\
(ii) $t=1$ if and only if $G \cong {\bf{Z}}_{p^{2}}$ or $G \cong
E_1$; \\
(iii)$t=2$ if and only if $G \cong
{\bf{Z}}_{p}\times{\bf{Z}}_{p^{2}}, G\cong D$ or $G\cong
{\bf{Z}}_{p} \times E_1$;\\
(iv)$t=3$ if and only if $G \cong {\bf{Z}}_{p^{3}},
G\cong{\bf{Z}}_{p} \times {\bf{Z}}_{p} \times{\bf{Z}}_{p^{2}},
G\cong D \times {\bf{Z}}_{p^{2}}, G\cong E_2, G\cong Q$ or $G\cong
{\bf{Z}}_{p} \times {\bf{Z}}_{p} \times E_1$.\\
Here ${\bf{Z}}_{p^{m}}$ denotes the cyclic group of order $p^{m}$,
$D$ denotes the dihedral group of order 8, $Q$ denotes the
quaternion group of order 8, $E_1$ denotes the extra special group
of order $p^{3}$ with odd exponent $p$, and $E_2$ denotes the
extra special group of order $p^{3}$ with odd exponent $p^{2}$.\\
\end{Theorem}

Now, in this paper, we are going to find some structures for the
$p$-group $G$ when $G$ is abelian and $|\Phi(G)|=p^{a}$ with
conditions on the exponents of $G, M(G),$ and $S_2M(G)$. The
following useful theorem of I. Schur is frequently used in our
method.
\begin{Theorem}(I. Schur [5]). Let $G \cong
{\bf{Z}}_{n_1}\oplus{\bf{Z}}_{n_2}\oplus...\oplus{\bf{Z}}_{n_k}$,
where $n_{i+1}|n_i$ for all $i \in {1,2,...,k-1}$ and $k\geq 2$,
and let ${\bf{Z}}_{n}^{(m)}$ denote the direct product of $m$
copies of ${\bf{Z}}_{n}$. Then $$M(G)\cong
{\bf{Z}}_{n_2}\oplus{\bf{Z}}_{n_3}^{(2)}\oplus...\oplus{\bf{Z}}_{n_k}^{(k-1)}$$\end{Theorem}

\textbf{Remark 1.3}. Let $G$ be an abelian group with a free
presentation $F/R$. Since $F' \leq R$, ${\mathcal
N}_{c_1}M(G)=\gamma _{c_1+1}(F)/[R,\ _{c_1}F]$. Now, we can
consider $\gamma _{c_1+1}(F)/[R,\ _{c_1}F]$ as a free
presentation for ${\mathcal N}_{c_1}M(G)$ and hence$${\mathcal
N}_{c_2}M({\mathcal N}_{c_1}M(G))= \frac{\gamma _{c_2+1}(\gamma
_{c_1+1}(F))}{[R,\ _{c_1}F,\ _{c_2}\gamma _{c_1+1}F]}.$$
Therefore by (1) we have
$${\mathcal N}_{c_1,c_2}M(G)={\mathcal N}_{c_2}M({\mathcal
N}_{c_1}M(G)).$$ By continuing the above process we can show
that$${\mathcal N}_{c_1,c_2 ...,c_t}M(G)={\mathcal
N}_{c_t}M(...{\mathcal N}_{c_2}M({\mathcal N}_{c_1}M(G))...).$$
In particular, if $c_1=c_2=1$, then we have $S_2M(G)=M(M(G)).$
\section{Main Results}

 Through out the paper we assume that $G$ is an abelian $p$-group
of order $p^{n}$ with $|M(G)|=p^{\frac{n(n-1)}{2}-t}.$
\begin{Lemma}. Let $\Phi(G)$, the Frattini
subgroup of $G$, be of order $p^{a}$. Then
$n=(a(a+1)+2t+2m)/2a$, for some $m \in {\bf{N}}_0$.\\
\end{Lemma}
\textit{Proof}. Let $G={\bf {Z}}_{p^{\alpha_1}}\oplus{\bf
{Z}}_{p^{\alpha_2}}\oplus ... \oplus {\bf
{Z}}_{p^{\alpha_{n-a}}},$ where $\alpha_1 \geq \alpha_2 \geq
...\geq \alpha_{n-a}$. By Theorem 1.2 , $M(G)\cong {\bf
{Z}}_{p^{\alpha_2}}\oplus{\bf {Z}}_{p^{\alpha_3}}^{(2)}\oplus ...
\oplus {\bf {Z}}_{p^{\alpha_{n-a}}}^{(n-a-1)}$ and so
$|M(G)|=p^{\alpha_2+2\alpha_3+...+(n-a-1)\alpha_{n-a}}$. But
$M(G)$ has order $p^{\frac{n(n-1)}{2}-t}$. Therefore\\
\begin{eqnarray*}
\frac{n(n-1)}{2}-t&=&
\alpha_2+2\alpha_3+...+(n-a-1)\alpha_{n-a}\\&\geq&
1+2+...+(n-a-1)\\&=& \frac{(n-a)(n-a-1)}{2}\\&=&
\frac{n^{2}-(2a+1)n+a(a+1)}{2}.
\end{eqnarray*}
Hence $2an \geq 2t+a(a+1)$, and the result holds.
\begin{Lemma}. With the assumption and
notation of the previous lemma we have the following inequalities
for the exponent of $G$, $$p^{a-m+1}\leq exp(G)\leq p^{a+1}.$$
\end{Lemma}
\textit{Proof.} Clearly $G/\Phi(G)$ is an elementary abelian
$p$-group of order $p^{n-a}$ and so $ exp(G)\leq
p^{a+1}$.\\

 For the other inequality let $G \cong
{\bf{Z}}_{p^{\alpha_1}}\oplus{\bf{Z}}_{p^{\alpha_2}}\oplus...
\oplus{\bf{Z}}_{p^{\alpha_{n-a}}}$,
where $\alpha_1\geq\alpha_2\geq...\geq\alpha_{n-a}\geq 1$. Then
similar
to the proof of previous lemma we have \\
\begin{eqnarray*}
\frac{n(n-1)}{2}-t&=&\alpha_2+2\alpha_3+...+(n-a-1)\alpha_{n-a}\\&
\geq& \alpha_2+(2+3+...+n-a-1) \\& \geq&
 \frac{2\alpha_2-2+n^{2}-(2a+1)n+a(a+1)}{2}\\
\end{eqnarray*}
Therefore $n\geq\frac{a(a+1)+2\alpha_2-2+2t}{2a}$ and hence by
Lemma 2.1 we have $\alpha_2\leq m+1$.\\Now, suppose by contrary
$exp(G)=p^{a-k+1}$, where $k>m$, then $\alpha_3\geq 2$. Thus by
Theorem 1.2 we have
\begin{eqnarray*}
\frac{n(n-1)}{2}-t&=&
\alpha_2+2\alpha_3+...+(n-a-1)\alpha_{n-a}\\&=&
(\alpha_2+\alpha_3+...+\alpha_{n-a})+\alpha_3+2\alpha_4+...+(n-a-2)\alpha_{n-a}
\\&\geq& (n-a+k-1)+(2+2+3+...+n-a-2)\\&=& n-a+k+\frac{(n-a-2)(n-a-1)}{2}. \\
\end{eqnarray*}
Hence $n\geq \frac{2k+a(a+1)+2+2t}{2a}>\frac{2m+a(a+1)+2+2t}{2a}$
which is a contradiction by Lemma 2.1.
\begin{Theorem}. With the above notation and
assumptions, let $G$ be of exponent $p^{a-m+1}$. Then
 $G\cong{\bf
{Z}}_{p^{a-m+1}}\oplus{\bf {Z}}_{p^{m+1}}\oplus \underbrace{{\bf
{Z}}_p\oplus...\oplus{\bf {Z}}_p}_{n-a-2-copies}.$\\
\end{Theorem}
\textit{Proof}. let $G \cong
{\bf{Z}}_{p^{\alpha_1}}\oplus{\bf{Z}}_{p^{\alpha_2}}\oplus...
\oplus{\bf{Z}}_{p^{\alpha_{n-a}}}$, where
$\alpha_1\geq\alpha_2\geq...\geq\alpha_{n-a}\geq 1$. By the proof
of previous lemma we have $\alpha_2 \leq m+1$. If $\alpha_2 \leq
m$, then $\alpha_3\geq2$ and we have
\begin{eqnarray*}
\frac{n(n-1)}{2}-t&=&
\alpha_2+2\alpha_3+...+(n-a-1)\alpha_{n-a}\\&=&
(\alpha_2+\alpha_3+...+\alpha_{n-a})+\alpha_3+2\alpha_4+...+(n-a-2)\alpha_{n-a}
\\&\geq& (n-a+m-1)+(2+2+3+...+n-a-2)\\&=& n-a+m+\frac{(n-a-2)(n-a-1)}{2}. \\
\end{eqnarray*}
Therefore $n\geq\frac{2m+a(a+1)+2+2t}{2a}$ which is a
contradiction by Lemma 2.1. Hence the result holds.
\begin{Theorem}. Further to the previous
notation and assumptions, let $m=k+s\ \  (k, s \in {\bf{N}}_0),$
$exp(G)=p^{a-k+1}$,
and $exp(M(G))+exp(S_2M(G))=p^{k+r}.$ Then\\
$$G\cong{\bf {Z}}_{p^{a-k+1}}\oplus{\bf {Z}}_{p^{x}} \oplus{\bf
{Z}}_{p^{k+r-x}} \oplus{\bf {Z}}_{p^{h_1}} \oplus...\oplus{\bf
{Z}}_{p^{h_f}}\oplus \underbrace{{\bf {Z}}_p \oplus...\oplus{\bf
{Z}}_p}_{n-a-(f+3)-copies},$$
where\\
$x=k-s+2r-3+3(h_1-1)+...+(f+2)(h_f-1)$, $h_1\geq h_2\geq ...\geq
h_f\geq 2$, and $f \leq -r+2$.\\
\end{Theorem}
\textit{Proof}.  Let $G \cong {\bf {Z}}_{p^{\alpha_1}}\oplus{\bf
{Z}}_{p^{\alpha_2}} \oplus...\oplus{\bf {Z}}_{p^{\alpha_{n-a}}},$
where $\alpha_1 \geq \alpha_2 \geq...\geq \alpha_{n-a}.$ By
Theorem 1.2 and Remark 1.3 it is easy to see that $exp(M(G))=
p^{\alpha_2}$, $exp (S_2M(G))=p^{\alpha_3}$, and so by hypothesis
$ \alpha_2+\alpha_3=k+r$, and $\alpha_1=a-k+1$. If $\alpha_3=1$,
then it is easy to see that $s=0$ and so $exp(G)=a-m+1.$ Hence by
Theorem 2.3 $G\cong{\bf {Z}}_{p^{a-m+1}}\oplus{\bf
{Z}}_{p^{m+1}}\oplus \underbrace{{\bf {Z}}_p\oplus...\oplus{\bf
{Z}}_p}_{n-a-2-copies}.$\\

Now, we can assume that $\alpha_3\geq2$ and hence $G$ may have the
following structure:\\
$$G\cong{\bf {Z}}_{p^{a-k+1}}\oplus{\bf {Z}}_{p^{x}} \oplus{\bf
{Z}}_{p^{k+r-x}} \oplus{\bf {Z}}_{p^{h_1}} \oplus...\oplus{\bf
{Z}}_{p^{h_f}}\oplus \underbrace{{\bf {Z}}_p \oplus...\oplus{\bf
{Z}}_p}_{n-a-(f+3)-copies},$$ where $f\geq0$ and $h_1\geq h_2\geq
...\geq h_f\geq 2$. If $f\geq1$ since $G$ is a group of order
$p^{n}$, we have\\
$(a-k+1)+(k+r)+(h_1+...+h_f)+(n-a-f-3)=n$ so $h_1+...+h_f=-r+f+2$
and hence $h_f=-r+f+2-h_1-...-h_{f-1}.$\\
But $h_1\geq h_2\geq ...\geq h_f\geq 2$, so that $-r+f+2\geq2f$
and so
$0\leq f\leq -r+2$.\\
Now, by Theorem 1.2 we have\\
$$\frac{n^{2}-n-2t}{2}=z+(1+2+...+n-a-1),$$
where$$z=x+2k+2r-2x+3h_1+4h_2+...+(f+2)h_{f}-(1+2+...+f+2).$$
Thus\\
$$n=\frac{2z+a(a+1)+2t}{2a}.$$ On the other hand by the
hypothesis $n=\frac{2(k+s)+a(a+1)+2t}{2a}$, hence we have $z=k+s$,
and the result follows.
\begin{Corollary}. With the notation and
assumptions of
previous theorem we have\\
(i) $G\cong{\bf {Z}}_{p^{a-k+1}}\oplus{\bf {Z}}_{p^{k-s+1}}
\oplus{\bf {Z}}_{p^{s+1}}\oplus \underbrace{{\bf
{Z}}_p\oplus...\oplus{\bf
{Z}}_p}_{n-a-3-copies},\ \ \ \ \ \ \ \ \ \  if \ \ \ r=2; $\\
(ii) $G\cong{\bf {Z}}_{p^{a-k+1}}\oplus{\bf {Z}}_{p^{k-s+2}}
\oplus{\bf {Z}}_{p^{s-1}} \oplus{\bf {Z}}_{p^{2}} \oplus
\underbrace{{\bf {Z}}_p \oplus...\oplus{\bf
{Z}}_p}_{n-a-4-copies},\ \ \ if\ \ \ r=1; $\\
(iii) if $r=0$, then \\
 $G \cong\left\{\begin{array}{ll}

 {\bf {Z}}_{p^{a-k+1}}\oplus{\bf {Z}}_{p^{k-s+4}} \oplus{\bf
{Z}}_{p^{s-4}} \oplus{\bf {Z}}_{p^{2}} \oplus{\bf
{Z}}_{p^{2}}\oplus \underbrace{{\bf {Z}}_p \oplus...\oplus{\bf
{Z}}_p}_{n-a-5-copies}&  \\

or\\

{\bf{Z}}_{p^{a-k+1}}\oplus{\bf {Z}}_{p^{k-s+3}} \oplus{\bf
{Z}}_{p^{s-3}} \oplus{\bf {Z}}_{p^{3}} \oplus \underbrace{{\bf
{Z}}_p \oplus...\oplus{\bf {Z}}_p}_{n-a-4-copies}& ;

\end{array}\right.$ \\
(iv) if $r=-1$, then\\
 $G \cong\left\{\begin{array}{ll}

 {\bf {Z}}_{p^{a-k+1}}\oplus{\bf {Z}}_{p^{k-s+8}} \oplus{\bf
{Z}}_{p^{s-9}} \oplus{\bf {Z}}_{p^{4}}\oplus \underbrace{{\bf
{Z}}_p \oplus...\oplus{\bf
{Z}}_p}_{n-a-4-copies}&  \\

{\bf{Z}}_{p^{a-k+1}}\oplus{\bf {Z}}_{p^{k-s+5}} \oplus{\bf
{Z}}_{p^{s-6}} \oplus{\bf {Z}}_{p^{3}} \oplus {\bf {Z}}_{p^{2}}
\oplus\underbrace{{\bf {Z}}_p \oplus...\oplus{\bf
{Z}}_p}_{n-a-5-copies}\\

{\bf{Z}}_{p^{a-k+1}}\oplus{\bf {Z}}_{p^{k-s+7}} \oplus{\bf
{Z}}_{p^{s-8}} \oplus{\bf {Z}}_{p^{2}} \oplus {\bf {Z}}_{p^{2}}
\oplus{\bf {Z}}_{p^{2}} \oplus \underbrace{{\bf {Z}}_p
\oplus...\oplus{\bf {Z}}_p}_{n-a-6-copies}&.
\end{array}\right.$
\end{Corollary}
\textit{Proof}. i) If $r=2$, then $f=0$. Therefore
$$G\cong{\bf {Z}}_{p^{a-k+1}}\oplus{\bf {Z}}_{p^{x}} \oplus{\bf
{Z}}_{p^{k+r-x}}\oplus \underbrace{{\bf {Z}}_p\oplus...\oplus{\bf
{Z}}_p}_{n-a-3-copies},$$ where $
x=k-s+1$. Hence the result follows.\\
ii) If $r=1$, then $f=0,1$. If $f=0$, then $n=a-k+1+k+1+n-a-3$
which is a contradiction. Then $f=1$ and
$$G\cong{\bf {Z}}_{p^{a-k+1}}\oplus{\bf {Z}}_{p^{x}} \oplus{\bf
{Z}}_{p^{k+r-x}} \oplus{\bf {Z}}_{p^{h_1}} \oplus \underbrace{{\bf
{Z}}_p \oplus...\oplus{\bf {Z}}_p}_{n-a-4-copies},$$ where $
h_1=-r+f+2=2$ and $x=k-s+2+6-(1+2+3)=k-s+2$. Hence the result
follows.\\
iii) If $r=0$, then $f=0,1,2$. If $f=0$, then $h_1+h_2+...+h_f=0$
but we have $h_1+h_2+...+h_f=-r+f+2=2$ which is a contradiction.
If $f=1$, then $x=k-s+4$ and if $f=2$, then $x=k-s+3$.\\
iv)By a routine calculation similar to (ii) the result holds.

Note that we can continue the above corollary for other integers
$r<-1$, but with a boring calculations.\\
 {\bf Acknowledgement}

 This research was in part supported by a grant from Center of
 Excellence in Analysis on Algebraic Structures, Ferdowsi University
 of Mashhad, CEAAS.

\end{document}